\documentclass[11pt]{article}
\usepackage[T1]{fontenc}
\usepackage{graphicx}
\usepackage{color}
\usepackage{amsmath,amssymb}
\usepackage{latexsym, amsthm}
\usepackage{enumerate}
\usepackage{bbm}
\usepackage{hyperref}
\usepackage{authblk}
\usepackage{setspace} 
\usepackage{accents}
\usepackage{ulem} 
\large\normalsize

\theoremstyle{plain}
\newtheorem{theorem}{Theorem}[section]

\newtheorem{lemma}[theorem]{Lemma}

\theoremstyle{definition}

\numberwithin{equation}{section}
\newcommand{\n}{\mathbb{N}}

\newcommand{\norma}[1]{\left\|#1\right\|}

\title{
The law of the iterated logarithm for a piecewise deterministic Markov process assured by the properties \\of the Markov chain given by the post-jump locations
}
\author[1]{Dawid Czapla}
\author[2]{Sander C. Hille}
\author[1]{Katarzyna Horbacz}
\author[1]{Hanna~Wojew\'odka-\'Sci\k{a}\.zko}
\affil[1]{\small Institute of Mathematics, University of Silesia in Katowice, Bankowa 14, 40-007 Katowice, Poland}
\affil[2]{\small Mathematical Institute, Leiden University, P.O. Box 9512, 2300 RA Leiden, The Netherlands}
\date{}
\begin{document}
\maketitle
\begin{abstract}
In the paper we consider some piecewise deterministic Markov process whose continuous component evolves according to semiflows, which are switched at the jump times of a Poisson process. The associated Markov chain describes the states of this process directly after the jumps. Certain ergodic properties of these two objects have already been investigated in our recent papers. We now aim to establish the law of the iterated logarithm for the aforementioned continuous-time process. Moreover, we intend to do this using the already proven properties of the discrete-time system. The abstract model under consideration has interesting interpretations in real-life sciences, such as biology. Among others, it can be used to describe the  stochastic dynamics of gene expression. 
\end{abstract}

\section*{Introduction}\label{sec:intro}

The law of the iterated logarithm (LIL) characterises essentially the maximal fluctuations around the mean of a stochastic process in discrete or continuous time. It is intimately related to the strong law of large numbers (SLLN) and the central limit theorem (CLT). The history of results on the LIL dates back to the work by Khinchin \cite{khinchin}, in the specific context of dyadic representations of numbers, and to the one by Kolmogorov \cite{kolm}, for general sequences of independent, non-necessarily identically distributed random variables that satisfy a particular `asymptotic boundedness' condition. Kolomogorov's results for identically distributed random variables with finite second moment were further generalised into the version of the LIL known as the Hartman-Wintner Theorem \cite{har_win}. See also \cite{strassen} and e.g. \cite{bingham} for a review of results on the LIL for the case of independent variables at the time of writing.

The main goal of this paper is to prove the validity of the LIL for a class of piecewise deterministic Markov processes (PDMPs). In this setting, the associated random variables are neither independent, nor identically distributed. Our method of proof is intentionally such that the result for the PDMP is derived from the validity of the LIL for the Markov chain given by its post-jump locations. The latter has been established in \cite{lil_chw} (see also the references mentioned there).

PDMPs have been introduced by Davis \cite{davis} as a general class of stochastic processes. They are encountered as suitable mathematical models for processes in the physical world around us, e.g. in  biology, as stochastic model for gene expression \cite{mtky}, gene regulation \cite{hhs}, excitable membranes \cite{riedler} or population dynamics \cite{alkurdi2, alkurdi}, as well as in resource allocation and service provisioning (queing, cf. \cite{davis}). Questions of ergodicity and asymptotic stability of PDMPs defined on locally compact state spaces have been studied in detail in \cite{benaim2, benaim,costa,costa2000}. The case of non-locally compact state space has been studied much less so far (see e.g.   \cite{alkurdi,dawid,hhs,riedler,hw}). A similar statement applies to the study of limit theorems (see \cite{zw,riedler}). For more information on the validity of limit theorems (SLLN, CLT or LIL) for non-stationary processes one may consult \cite{bms,dawid,clt_chw,lil_chw,klo}.

A PDMP consists of deterministic movement in a state space (a~Polish metric space in our case) that is alternated at random times of intervention with a~random jump in state. In general, the distribution of the next intervention time and the jump can be both state dependent (cf. e.g. \cite{hhs}). Here, and e.g. also in \cite{alkurdi}, only the jump is distributed conditionally given the current state of the system.  
The process examined in this paper (described in detail in Section \ref{sec:}) involves jumps that occure at random time points according to a~Poisson process. Any post-jump location is attained by transforming a pre-jump state using a randomly selected function, and, further, by adding a radom shift to the resulting state. Between any two consecutive jumps, the system is driven deterministically by one of a finite number of flows, which are switched at the jump times.  
If the state space is augmented with an index set of the applied movements, then the chain obtained by pairing the state just after the jump with the index of the movement that is applied thereafter yields a Markov chain, which intuitively should contain `all information' about the PDMP. Therefore it is enlightning to show how properties of the PDMP can be proven from relevant properties of the Markov chain constituted by the post-jump states.

Essentially, our method of proof splits the problem into subproblems that can be analyzed separately. One subproblem can be addressed using a version of the LIL for certain square integrable martingales, whose proof draws heavily on \cite[Theorem 1]{hs} and uses the coupling methods applied for establishing \cite[Lemma 2.2]{clt_chw} (cf. also \cite{hairer}). Another builds on the validity of the LIL for Markov chains associated to PDMPs in the abstract model class, that has been obtained recently (cf. \cite[Theorem 4.1]{lil_chw}).

We believe that the class of dynamical systems under study is broad enough to cover models of suitable real-life systems, e.g. biological systems, such as artificial evolutionary experiments on bacteria \cite{alkurdi2}, as well as chemotactic movement of bacteria or amoeba (related to the study of so-called velocity-jump models, employing particular Fokker-Planck equations, see e.g. \cite{oda, ho, perthame}). Discussion and the detailed study of such application are beyond the scope of this paper, but they shall be the subject of our further reseach collaboration.

\section{Prelimenaries}\label{sec:1}

Let us first introduce a piece of notation, as well as gather the most important definitions and facts, used in this paper.

\subsection{Some notation and basic definitions}
For any point $x$ and any set $A$, the symbols $\delta_x$ and $\mathbbm{1}_A$ will denote the Dirac measure at $x$~and the indicator function of $A$, respectively.

Suppose that $(E,\varrho_E)$ is a Polish metric space and let $\mathcal{B}_E$ denote the $\sigma$-field of its Borel subsets. Let $B_b(E)$ stand for the space of all bounded, Borel measurable functions \hbox{$f:E\to\mathbb{R}$} equipped with the supremum norm $\|f\|_{\infty}=\sup_{x\in E}|f(x)|$. We shall also refer to certain subspaces of $B_b(E)$, namely $C_b(E)$, consisting of all continuous functions, $Lip_b(E)$, consisting of all Lipschitz continuous functions,
and 
$$Lip_{FM}(E)=\{f\in Lip_b(E):\,\|f\|_{BL}\leq 1\},$$ where the norm $\|\cdot\|_{BL}$ is given by \hbox{$\|f\|_{BL}=\max\{|f|_{Lip},\,\norma{f}_{\infty}\}$} and $|f|_{Lip}$ stands for the minimal Lipschitz constant of~$f$ for every $f\in Lip_b(E)$. 
Finally, we will also consider the space $\bar{B}_b(E)$ of functions $f:E\to\mathbb{R}$ which are Borel measurable and bounded below.

The spaces of finite and probability Borel measures on $E$ will be denoted by $\mathcal{M}_{fin}(E)$ and $\mathcal{M}_{1}(E)$, respectively. Further, we also define   
$$\mathcal{M}_{1,r}^V(E)=\left\{\mu\in \mathcal{M}_1(E):\, \int_EV^r(x)\,\mu(dx)<\infty\right\}$$ for any $r>0$ and any given Lyapunov function $V:E\to[0,\infty)$, that is, a function which is continuous, bounded on bounded sets, and, in the case of unbounded $E$, satisfies \linebreak\hbox{$\lim_{\varrho_E(x,\bar{x})\to\infty}V(x)=\infty$} for some fixed point $\bar{x}\in E$. 
For brevity, for any $f\in \bar{B}_b(E)$ and any signed Borel measure $\mu$ on $E$, we will write $\langle f,\mu\rangle$ for $\int_Ef(x)\,\mu(dx)$. As usual, $\text{supp}\,\mu$ will stand for the support of $\mu\in\mathcal{M}_{fin}(E)$.

To evaluate the distance between probability measures, we will use the so-called Fortet-Mourier distance (see e.g. \cite{l_frac}), defined as follows:
\begin{align*}
d_{FM}(\mu_1,\mu_2)=\sup\left\{\left|\left\langle f,\mu_1-\mu_2\right\rangle\right|:\;f\in Lip_{FM}(E)\right\}\quad\text{for}\quad \mu_1,\mu_2\in \mathcal{M}_1(E).
\end{align*}
Let us indicate that, under the assumption that $(E,\varrho_E)$ is a Polish space, the convergence in $d_{FM}$ is equivalent to the weak convergence of probability measures, and also the space $(\mathcal{M}_1(E),d_{FM})$ is complete (for the proofs of both these facts see e.g. \cite{dudley}).

\subsection{Markov operators and the semigroups of Markov operators}
A function 
$P:E\times \mathcal{B}_E\to [0,1]$ is called a (sub)stochastic kernel, if, for any fixed $A\in \mathcal{B}_E$, $P(\cdot,A):E\to[0,1]$ is a Borel measurable map, and, for any fixed $x\in E$, \hbox{$P(x,\cdot):\mathcal{B}_E\to[0,1]$} is a (sub)probability Borel measure. 
For any two kernels $P:E\times\mathcal{B}_E\to[0,1]$ and \linebreak\hbox{$R:E\times \mathcal{B}_E\to[0,1]$} we can define their composition $PR:E\times\mathcal{B}_E\to[0,1]$ given by 
\begin{align}\label{eq:composition}
PR(x,A)=\int_EP(y,A)R(x,dy)\;\;\;\text{for}\;\;\;x\in E\;\;\;\text{and}\;\;\;A\in\mathcal{B}_E.
\end{align}
Following this rule, for any (sub)stochastic kernel $P:E\times\mathcal{B}_E\to[0,1]$, we can define its $n$-th step kernels $P^n:E\times\mathcal{B}_E\to[0,1]$, inductively on $n\in\mathbb{N}$, by setting
$P^n=PP^{n-1}$, where $P^0$ is given by $P^0(x,A)=\delta_x(A)$ for every $x\in E$ and any $A\in\mathcal{B}_E$. 

Moreover, for any stochastic kernel $P$, we can define a regular Markov operator \linebreak\hbox{$(\cdot)P:\mathcal{M}_{fin}(E)\to \mathcal{M}_{fin}(E)$} and its dual operator \hbox{$P(\cdot):{B}_b(E)\to {B}_b(E)$} in the following way:
\begin{align}\label{def:markov_op}
&\mu P(A)=\int_EP(x,A)\,\mu(dx)\;\;\;\text {for} \;\;\;\mu\in \mathcal{M}_{fin}(E),\; A\in \mathcal{B}_E,\\
\label{def:dual_op}
&Pf(x)=\int_Ef(y)\,P(x,dy) \;\;\;\text{for}\;\;\; f\in {B}_b(E),\;x\in E.
\end{align}
Obviously, $\langle f,\mu P\rangle=\langle Pf,\mu\rangle$ for any $f\in{B}_b(E)$ and any $\mu\in \mathcal{M}_{fin}(E)$. Moreover, note that any operator $P(\cdot)$ of the form \eqref{def:dual_op} can be extended, in the usual way, to a linear operator on $\bar{B}_b(E)$, preserving the duality property, and hence it is reasonable to apply $P(\cdot)$ to any Lyapunov function. For notational simplicity, we shall use the same symbol for the extension as for the original operator on ${B}_b(E)$. 
An operator $(\cdot)P$, given by \eqref{def:markov_op}, is said to be Markov-Feller if $Pf\in C_b(E)$ for every $f\in C_b(E)$. 

We call $\mu_*\in \mathcal{M}_{fin}(E)$ an~invariant measure of $(\cdot)P$ if $\mu_* P =\mu_*$.  If~
$(\cdot)P$ has a unique invariant measure $\mu_*\in \mathcal{M}_1(E)$ and there exists $q\in (0,1)$ such that
\begin{align*}
d_{FM}(\mu P^n,\mu_*)\leq c(\mu)q^n \;\;\;\mbox{for any}\;\;\;\mu\in\mathcal{M}_{1,1}^V(E),\;n\in\n,
\end{align*}
where $c(\mu)$ is a constant depending only on $\mu$, then $(\cdot)P$ is said to be exponentially ergodic in~$d_{FM}$.

Let us consider $E^{\mathbb{N}_0}$ with the product topology. For every $n\in\mathbb{N}_0$ define $\phi_n:E^{\mathbb{N}_0} \to E$ by the formula $\phi_n(\omega)=e_n$, where $\omega=(e_0,e_1,\ldots)\in E^{\mathbb{N}_0}$. According to \hbox{\cite[Theorem 2.8]{revuz}}, for any \hbox{$\mu\in \mathcal{M}_{1}(E)$} and any stochastic kernel~$P:E\times \mathcal{B}_E\to[0,1]$, there exists \hbox{$\mathbb{P}\in\mathcal{M}_1(E^{\mathbb{N}_0})$} such that $(\phi_n)_{n\in\n_0}$ is a time-homogeneus Markov chain on the probability space $(E^{\mathbb{N}_0},\mathcal{B}_{E^{\mathbb{N}_0}},\mathbb{P})$ with transition function $P$ and initial measure $\mu$, that is 
\begin{align}\label{eq:star}
P^n(x,A)=\mathbb{P}(\phi_{k+n}\in A | \phi_k=x)\;\;\;\text{for}\;\;\; x\in E,\; A\in\mathcal{B}_E,\; n,k\in\mathbb{N}_0,
\end{align}
and
\[\mu(A) = \mathbb{P}(\phi_0\in A) \;\;\;\text{for}\;\;\;A\in\mathcal{B}_E.\]
The chain defined as above shall be further called the canonical Markov chain. Clearly, 
$\mathbb{P}(B)$ may be read as the probability of the event $\left\{(\phi_n)_{n\in\mathbb{N}_0}\in B\right\}$ for any $B\in\mathcal{B}_{E^{\mathbb{N}_0}}$.

Conversely, it is clear that the one-step transition law of any time-homogeneous Markov chain determines a stochastic kernel and the corresponding $n$-step kernels which satisfy \eqref{eq:star}.

As far as the dual operator $P(\cdot)$ is concerned, we have 
\[P^nf(x)=\mathbb{E}\left(f(\phi_n)|\phi_0=x\right)\;\;\;\text{for}\;\;\;x\in E,\;f\in B_b(E),\;n\in\mathbb{N}.\]

A regular Markov semigroup $({P}_t)_{t\in\mathbb{R}_+}$ is a family of regular Markov operators \linebreak \hbox{$(\cdot){P}_t: \mathcal{M}_{fin}(E) \to \mathcal{M}_{fin}(E)$}, $t\in\mathbb{R}_+$, which form a semigroup (under composition) with the identity transformation $(\cdot){P}_0$ as the unity element. 
Provided that $(\cdot){P}_t$ is a Markov-Feller operator for every $t\in\mathbb{R}_+$, the semigroup $({P}_t)_{t\in\mathbb{R}_+}$ is said to be Markov-Feller, too. If, for some $\mu_*\in\mathcal{M}_{fin}(E)$, $\mu_*{P}_t=\mu_*$ for every $t\in\mathbb{R}_+$, then we call $\mu_*$ an invariant measure of $({P}_t)_{t\in\mathbb{R}_+}$.

Let $(\phi(t))_{t\in\mathbb{R}_+}$ be an $E$-valued time-homogeneous Markov process, defined on an arbitrary probability space $(\Omega,\mathcal{F},\mathbb{P})$, with continuous time parameter $t\in\mathbb{R}_+$. 
Suppose that, for any $t\in\mathbb{R}_+$,  \hbox{$P_t:E\times\mathcal{B}_E\to[0,1]$} is defined by 
\begin{equation}\label{trans_t} 
P_t(x,A)=\mathbb{P}(\phi(t)\in A| \phi(0)=x)\;\;\;\mbox{for}\;\;\;x\in E,\;A\in\mathcal{B}_E,\;t\in\mathbb{R}_+.
\end{equation}
It is well-known that these transition probability functions form a semigroup of stochastic kernels under the composition operation defined by \eqref{eq:composition}. Thus the family of the corresponding Markov operators $(P_t)_{t\in\mathbb{R}_+}$ is a regular Markov semigroup. 
The dual operator of $P_t$, $t\in\mathbb{R}_+$, can be expressed in the form 
\[P_tf(x)={\mathbb{E}}\left(f(\phi(t))|\phi(0)=x\right).\]

Now, let $(\phi_n)_{n\in\mathbb{N}_0}$ be a Markov chain with transition law $P$, and let $(\phi^{(1)}_n)_{n\in\mathbb{N}_0}$, $(\phi^{(2)}_n)_{n\in\mathbb{N}_0}$ be its copies with initial distributions $\mu_1\in\mathcal{M}_1(E)$, $\mu_2\in\mathcal{M}_1(E)$, respectively. 
A~time-homogeneus Markov chain $(\phi^{(1)}_n,\phi^{(2)}_n)_{n\in\mathbb{N}_0}$ evolving on $E^2$ (endowed with the product topology) is said to be a Markovian coupling of $(\phi^{(1)}_n)_{n\in\mathbb{N}_0}$ and $(\phi^{(2)}_n)_{n\in\mathbb{N}_0}$ whenever its transition law $C:E^2\times  \mathcal{B}_{E^2}\to\left[0,1\right]$ satisfies \[C((x,y),A\times E)= P(x,A)\;\;\;\text{and}\;\;\;C((x,y),E\times A)= P(y,A)\;\;\;\text{for any}\;\;\;x,y\in E,\; A\in \mathcal{B}_E,\]
and its initial distribution $\alpha \in \mathcal{M}_1(E^2)$ is such that
\[\alpha(A \times E)=\mu_1(A), \;\;\;\alpha (E \times A)=\mu_2(A)\;\;\;\text{for any}\;\;\;A\in\mathcal{B}_E.\] 
In what follows we always assume that the coupling is defined canonically on the coordinate space $((E^2)^{\mathbb{N}_0},\mathcal{B}_{(E^2)^{\mathbb{N}_0}})$ endowed with an appropriately constructed measure $\mathbb{C}\in\mathcal{M}_1((E^2)^{\mathbb{N}_0})$.

\subsection{The law of the iterated logarithm for Markov processes}\label{sec:def_lil}

Consider an $E$-valued time-homogeneous Markov chain $(\phi_n)_{n\in\mathbb{N}_0}$ with initial distribution $\mu\in\mathcal{M}_1(E)$ and an $E$-valued time-homogeneous Markov process $({\phi}(t))_{t\in\mathbb{R}_+}$ with initial distribution $\nu\in\mathcal{M}_1(E)$.  
For any function $g\in Lip_b(E)$, let us introduce \hbox{$(s_n({g}))_{n\in\mathbb{N}_0}$} and \hbox{$(s(g)(t))_{t\in\mathbb{R}_+}$}, given by 
\begin{gather}\label{def:s_n}
s_n({g})=\cfrac{\sum_{i=0}^{n-1}{g}(\phi_i)}{\sqrt{2n\ln(\ln (n))}} 
\;\;\;\text{for} \;\;\;n>e\qquad\text{and}\qquad
s_n({g})=0\;\;\;\text{for} \;\;\;n\leq e;\\
\label{def:s_t}
{s}({g})(t)=\cfrac{\int_0^t{g}({\phi}(s))ds}{\sqrt{2t\ln(\ln (t))}} 
\;\;\;\text{for} \;\;\;t>e\qquad\text{and}\qquad
s(g)(t)=0\;\;\;\text{for} \;\;\;t\leq e.
\end{gather}
Suppose that \hbox{$\mu_*\in\mathcal{M}_1(E)$} and $\nu_*\in\mathcal{M}_1(E)$ are the unique invariant measures for $(\phi_n)_{n\in\mathbb{N}_0}$ and $({\phi}(t))_{t\in\mathbb{R}_+}$, respectively. 
We say that the Markov chain $\left(g(\phi_n)\right)_{n\in\mathbb{N}_0}$ satisfies the LIL if, for $\hat{g}=g-\langle g,\mu_*\rangle$ and some $\sigma(\hat{g})\in(0,\infty)$,  
\[\limsup\limits_{n\to\infty}s_n\left(\hat{g}\right)=\sigma(\hat{g})
\;\;\;
\text{and}\;\;\;
\liminf\limits_{n\to\infty}s_n\left(\hat{g}\right)=-\sigma(\hat{g})\;\;\;\mathbb{P}\text{-a.s.}\]
Accordingly, we say that the Markov process $\left(g({\phi}(t))\right)_{t\in\mathbb{R}_+}$ satisfies the LIL if, for \linebreak\hbox{$\bar{g}=g-\langle g,\nu_*\rangle$} and some $\overline{\sigma}(\bar{g})\in(0,\infty)$,
\[\limsup\limits_{t\to\infty}s\left(\bar{g}\right)(t)=\overline{\sigma}(\bar{g})
\;\;\;
\text{and}\;\;\;
\liminf\limits_{t\to\infty}s\left(\bar{g}\right)(t)=-\overline{\sigma}(\bar{g})\;\;\;{\mathbb{P}}\text{-a.s.}\]

\section{An abstract model}\label{sec:}

In the beginning, we shall discuss the structure and assumptions of the model under consideration. Let us indicate that this model was initially introduced in \cite{dawid}, where we have also elaborated on its possible applications. Further, let us summarise the already known results that are used further in this paper.

\subsection{The structure of the model and the undertaken assumptions}

Consider a separable Banach space $(H,\|\cdot\|)$ and a closed subset $Y$ of $H$. For any $x\in H$ and any $r>0$, let $B(x,r)$ denote an open ball in $H$ centered at $x$ and of radius $r$. Let us also fix a~topological measure space $(\Theta,\mathcal{B}_{\Theta},\vartheta)$ with a finite Borel measure~$\vartheta$. With a slight abuse of notation, we will further write $d\theta$ only, instead of $\vartheta(d\theta)$. Finally, fix $m\in\n$ and introduce the set of indexes $I:=\{1,\ldots,m\}$ equipped with the metric $d$ given by 
\begin{align*}
d(i,j)=\left\{
\begin{array}{ll}
1,&i\neq j\\
0,&i=j
\end{array}.
\right.
\end{align*}

We shall investigate a random dynamical system $({Y}(t))_{t\in\mathbb{R}_+}$ evolving through jumps, occuring at random moments $\tau_n$, $n\in\mathbb{N}$, which coincide with the jump times of a Poisson process with a given intensity $\lambda$. 
In every time interval $[\tau_{n-1},\tau_{n})$, where $\tau_0=0$, the system is driven by one of the given continuous semiflows \hbox{$S_i:\mathbb{R}_+\times Y\to Y$}, $i\in I$. 
The current semiflow, say $S_i$, is switched to another (or the same) one $S_j$ with a probability $\pi_{ij}(y)$, depending on the post-jump state $y$. We assume that these place-dependent probabilities constitute a matrix of continuous functions  $\pi_{ij}: Y\to [0,1]$, $i,j\in I$, such that
\[\sum_{j\in I}\pi_{ij}(y)=1\;\;\;\text{for any}\;\;\;y\in Y,\,i\in I.\]
The above description can be shortly formalized by the following formula:
\begin{align}\label{def:Y(t)}
{Y}(t)=S_{\xi_n}\left(t-\tau_n,Y_n\right)\;\;\; \text{for}\;\;\; t\in[\tau_n,\tau_{n+1}),
\end{align}
where $\xi_n$ is an $I$-valued random variable  indicating which semiflow has been chosen after the $n$-th jump, and $Y_n$ is a result of some transformation of the state $Y(\tau_n-)$ just before the jump. The transformation is attained by a function \hbox{$w_{\theta}:Y\to Y$}, selected randomly among all possible ones $\{w_{\theta}: \theta\in \Theta\}$, and further disturbed by adding some random shift $H_n$. Therefore, we can formally write 
\begin{align*}
Y_n=w_{\theta_n}({Y}(\tau_n-))+H_n.
\end{align*}

It is assumed that, given ${Y}(\tau_n-)=y$, the probability of choosing $w_{\theta}$ (at the jump time $\tau_n$) is determined by the density function $\theta \mapsto p(y,\theta)$ such that \hbox{$p: Y\times \Theta \to \left[0,\infty\right)$} is a~continuous map. 
Moreover, it is required that the map $(y,\theta)\mapsto w_{\theta}(y)$ is continuous. Further, we also assume that, for some $\varepsilon>0$, all the variables $H_n$, $n\in\mathbb{N}$,  have a common distribution \hbox{$\nu^{\varepsilon}\in\mathcal{M}_1(H)$} supported on $B(0,\varepsilon)\subset H$, and that  
\[w_{\theta}(y)+h\in Y\;\;\; \text{for any} \;\;\; h\in\text{supp}(\nu^{\varepsilon}),\; \theta \in \Theta,\;y\in Y.\]

We therefore formally consider a stochastic process 
$({Y}(t))_{t\in\mathbb{R}_+}$ of the form \eqref{def:Y(t)}, defined as an interpolation of the discrete-time process $(Y_n)_{n\in\mathbb{N}_0}$ determined  
by the recursive formula
\begin{equation}\label{def:Y_n}
Y_{n}=Y(\tau_n)=w_{\theta_{n}}(S_{\xi_{n-1}}(\Delta \tau_{n},Y_{n-1}))+H_{n}\;\;\;\text{for}\;\;\;n\in\mathbb{N},
\end{equation}
where $(\tau_n)_{n\in\n_0}$, $(\theta_n)_{n\in\n}$, $(\xi_n)_{n\in\n_0}$ and $(H_n)_{n\in\n}$ are certain sequences of random variables (specified below) with values in $\mathbb{R}_+$, $\Theta$, $I$ and $H$, respectively.

The distribution of $(Y_0,\xi_0)$ is fixed arbitrarily. The sequence $(\tau_n)_{n\in\mathbb{N}_0}$, wherein $\tau_0=0$ a.s., is such that $\tau_n\to\infty$ a.s., as $n\to\infty$. The increments $\Delta \tau_{n+1}:=\tau_{n+1}-\tau_{n}$, $n\in\mathbb{N}_0$, are, in turn, assumed to be mutually independent and identically distributed according to the exponential distribution with intensity $\lambda>0$. 
Moreover, the disturbances $(H_n)_{n\in\n}$ are required to be identically distributed with~$\nu^{\varepsilon}$, introduced above.  
Finally, the chains $(\xi_n)_{n\in\n}$ and $(\theta_n)_{n\in\n}$ are defined, inductively on $n\in\mathbb{N}_0$, as follows:
\begin{align*}
&\mathbb{P}(\xi_{n+1}=j\;|\;Y_{n+1}=y,\, \xi_{n}=i;\,W_n)=\pi_{ij}(y) \;\;\;\text{for}\;\;\; y\in Y,\;i, j\in I,\\
&\mathbb{P}(\theta_{n+1}\in D\;|\;S_{\xi_{n}}(\Delta \tau_{n+1},Y_{n})=y;\,W_n)
=\int_{D} p(y,\theta)\,d\theta\;\;\;\text{for}\;\;\; D\in\mathcal{B}_{\Theta},\;y\in Y, 
\end{align*}
where
$$W_0=(Y_0,\;\xi_0)\;\;\,\text{and}\;\;\,W_n=(W_0,\;H_1,\ldots,H_n,\;\tau_1,\ldots,\tau_n,\;\theta_1,\ldots,\theta_n,\;\xi_1,\ldots,\xi_n)\;\;\,\text{for}\;\;\,n\in\mathbb{N}.$$
We also demand that, for any $n\in\mathbb{N}_0$, the variables $\Delta\tau_{n+1}$, $H_{n+1}$, $\theta_{n+1}$ and $\xi_{n+1}$ are ({mutually}) conditionally independent given $W_n$, and that $\Delta\tau_{n+1}$ and $H_{n+1}$ are independent of $W_n$.

Let us now consider the space $X:=Y\times I$ with the metric \hbox{$\varrho_{c}$}, 
given by
\begin{align}\label{def:rho_c}
\varrho_{c}\left((y_1,i_1),(y_2,i_2)\right)=\|y_1-y_2\|+c\,{d}(i_1,i_2)\quad\text{for}\quad (y_1,i_1),(y_2,i_2)\in X,
\end{align}
with a sufficiently large constant $c\geq 1$ (defined explictly in \cite{dawid}).
Now, define $$X_n:=(Y_n,\xi_n)\;\;\; \text{for}\;\;\;{n\in\mathbb{N}_0}.$$ 
Given $\mu\in \mathcal{M}_1(X)$, we shall further consider  the canonical $(X\times\mathbb{R}_+)$-valued Markov chain $(X_n,\Delta\tau_n)_{n\in\n_0}$ with initial distribution $\mu\otimes \delta_0$, defined on a probability space $(\Omega, \mathcal{F},\mathbb{P})$, where $\Omega:=(X\times\mathbb{R}_+)^{\mathbb{N}_0}$ and $\mathcal{F}:=\mathcal{B}_{\Omega}$, whose transition law $\Pi:(X\times\mathbb{R}_+)\times\mathcal{B}_{X\times\mathbb{R}_+}\to[0,1]$ is given by 
\begin{align}\label{def:Pi_transition}
\Pi\left((y,i,s),A\right)
=&\int_0^{\infty}\lambda e^{-\lambda t}\int_{\Theta}p(S_i(t,y),\theta) \int_{\text{supp}(\nu^{\varepsilon})}\Bigg(\sum_{j\in I}\mathbbm{1}_A(w_{\theta}(S_i(t,y))+h,j,t)
\nonumber\\
&\times \pi_{ij}(w_{\theta}(S_i(t,y))+h)\Bigg) \nu^{\varepsilon}(dh)\,d\theta\,dt
\end{align}
for any $(y,i,s)\in X\times\mathbb{R}_+$ and any $A\in\mathcal{B}_{X\times\mathbb{R}_+}$.  
Note that $(X_n)_{n\in\n_0}$ itself is also a time-homogeneous Markov chain with transition law \hbox{$P:X\times\mathcal{B}_X\to\left[0,1\right]$} satisfying
\begin{align}\label{def:Pi_epsilon}
P((y,i),A)=\Pi\left((x,s),A\times\mathbb{R}_+\right)\;\;\;\text{for any}\;\;\;x\in X,\;s\in\mathbb{R}_+\;\;\;\text{and}\;\;\;A\in\mathcal{B}_X.
\end{align}
Moreover, we have
\begin{align*}
\Pi\left((x,s),X\times B\right)=\int_B\lambda e^{-\lambda t}dt\;\;\;\text{for any}\;\;\;(x,s)\in X\times\mathbb{R}_+\;\;\;\text{and}\;\;\;B\in\mathcal{B}_{\mathbb{R}_+}.
\end{align*}

Now, define the continuous-time process $(X(t))_{t\in\mathbb{R}_+}$ on the space $(\Omega, \mathcal{F},\mathbb{P})$, by setting
\begin{align}\label{def:X(t)}
X(t)=\left(Y(t),\xi(t)\right)=\left(S_{\xi_n}\left(t-\tau_n,Y_n\right),\xi_n\right)\;\;\;\text{for}\;\;\;t\in[\tau_n,\tau_{n+1}).
\end{align}
One may check that $(X(t))_{t\in\mathbb{R}_+}$ is an $X$-valued time-homogeneous Markov process such that \[{X}(\tau_n)=X_n\;\;\;\text{for any}\;\;\;n\in\n_0.\] 
The Markov transition semigroup associated with the process $\left(X(t)\right)_{t\in\mathbb{R}_+}$ shall be denoted by $({P}_t)_{t\in\mathbb{R}_+}$.

Summarising this part of the paper, let us indicate that, if $X_0$ is distributed according to some measure \hbox{$\mu\in\mathcal{M}_1(X)$}, then we get
\begin{gather}\label{properties_P}
\begin{aligned}
\mathbb{P}\left(\left(X_n,\Delta\tau_n\right)\in D\right)=\left(\mu\otimes \delta_0\right)P^n(D)\;\;\;\text{for any}\;\;\;D\in\mathcal{B}_{X\times\mathbb{R}_+},
\\
\mathbb{P}\left(\Delta\tau_n\in B\right)=\int_B\lambda e^{-\lambda t}dt\;\;\;\text{for any}\;\;\;B\in\mathcal{B}_{\mathbb{R}_+},\;n\in\mathbb{N},\end{aligned}
\\ \label{properties_Pt}
\mathbb{P}\left(X(t)\in A\right)=\mu P_t(A)\;\;\;\text{for any}\;\;\;A\in\mathcal{B}_{X},\;t\in\mathbb{R}_+.
\end{gather}

Let us further assume that there exist $\bar{y}\in Y$, $\alpha\in\mathbb{R}$ and positive constants $L$, $\bar{L}$, $L_w$, $L_{\pi}$, $L_p$, $\delta_{\pi}$, $\delta_p$, $r\in(0,2)$ such that
\begin{equation}\label{eq:balance1} 
L^{2+r}L_w+(2+r)\frac{\alpha}{\lambda}<1, 
\end{equation}
and, for all $i,i_1,i_2\in I$, $y_1,y_2\in Y$, $t\in\mathbb{R}_+$, the following conditions hold:
\begin{gather}
\tag{A1}\label{cnd:a1}
\sup_{y\in Y}\int_0^{\infty}e^{-\lambda t}\int_{\Theta}\norma{ w_{\theta}(S_i(t,\bar{y}))-\bar{y}}^{2+r}p(S_i(t,y),\theta)\,d\theta\,dt<\infty,\\
\tag{A2}\label{cnd:a2}
\norma{S_{i_1}(t,y_1)-S_{i_2}(t,y_2)}\leq Le^{\alpha t}\norma{y_1-y_2}+t\bar{L}\,d(i_1,i_2),\\\tag{A3}\label{cnd:a3}
\int_{\Theta} p(y_1,\theta)\norma{w_{\theta}(y_1)-w_{\theta}(y_2)}^{2+r}\,d\theta\leq L_w\norma{y_1-y_2}^{2+r},\\
\tag{A4}\label{cnd:a4}
\sum_{j\in I} |\pi_{ij}(y_1)-\pi_{ij}(y_2)|\leq L_{\pi} \norma{y_1-y_2},\;\;\;\int_{\Theta} |p(y_1,\theta)-p(y_2,\theta)|\,d\theta\leq L_{p} \norma{y_1-y_2},\\
\tag{A5}\label{cnd:a5}
\sum_{j\in I} \min\{\pi_{i_1,j}(y_1),\pi_{i_2,j}(y_2)\}\geq \delta_{\pi},\;\;\;\int_{\Theta(y_1,y_2)}\min\{p(y_1,\theta),p(y_2,\theta)\}\,d\theta\geq \delta_p,
\end{gather}
where $\Theta(y_1,y_2):=\{\theta\in\Theta:\, \norma{w_{\theta}(y_1)-w_{\theta}(y_2)}\leq L_w \norma{y_1-y_2}\}$. 
Hypotheses \eqref{cnd:a1}-\eqref{cnd:a5} and their reasonableness are discussed in detail e.g. in \cite{dawid, lil_chw, asia}.

\subsection{Certain properties of the model under consideration}

Suppose that hypothesis \eqref{cnd:a1}-\eqref{cnd:a5} hold with constants satisfying \eqref{eq:balance1}. Then \linebreak\hbox{\cite[Theorem 4.1]{dawid}} implies that the Markov operator $P$, determined by \eqref{def:Pi_epsilon}, is exponentially ergodic in $d_{FM}$ induced by the metric \hbox{$\varrho_{c}$} given by \eqref{def:rho_c}. 
In fact, the exponential ergodicity itself can be obtained even under slightly weaker assumptions than \eqref{cnd:a1}-\eqref{cnd:a5} (cf. \cite{dawid}). To be more precise, \eqref{cnd:a1}, \eqref{cnd:a3} and  \eqref{eq:balance1} may be considered in their weaker versions, wherein $r=-1$. However, to establish the law of the iterated logarithm, we need them as given in \cite{lil_chw} and also in this paper.

Fix an arbitrary non-constant function $g\in Lip_b(X)$. Further, consider the chain $(X_n)_{n\in\mathbb{N}_0}$ governed by $P$, defined in \eqref{def:Pi_epsilon}, with the initial distribution $\mu\in\mathcal{M}_{1,2+r}^V(X)$, where $r\in(0,2)$ is the constant appearing in \eqref{eq:balance1}, and $V:X\to[0,\infty)$ is a Lyapunov function given by
\begin{align}
\label{def:V}
V(y,i)=\|y-\bar{y}\|\quad\text{for every}\quad(y,i)\in X,
\end{align}
where $\bar{y}$ is determined by \eqref{cnd:a1}.  
Referring to \cite[Theorem 4.1]{lil_chw}, we know that the chain $(g(X_n))_{n\in\mathbb{N}_0}$ satisfies the invariance principle for the LIL, and whence it also satisfies the LIL itself (cf. \cite[Section 3.2]{lil_chw}).

In \hbox{\cite[Corollary 4.5]{dawid}} we have proven that there is a one-to-one correspondence between invariant measures of the operator $P$ and those of the semigroup $({P}_t)_{t\in\mathbb{R}_+}$. This obviously implies that $({P}_t)_{t\in\mathbb{R}_+}$ has a unique invariant distribution if and only if $P$ admits the one, which holds, in particular, whenever conditions \eqref{cnd:a1}-\eqref{cnd:a5} and \eqref{eq:balance1} are satisfied. 
The above-mentioned correspondence can be described explicitly, using the Markov operators associated with the stochastic kernels $G,W: X\times \mathcal{B}_X\to\left[0,1\right]$ defined as follows:
\begin{gather}
\label{def_G}
G((y,i),A)=\int_0^{\infty} \lambda e^{-\lambda t} \mathbbm{1}_A(S_i(t,y),i)\,dt,\\
\label{def_W}
W ((y,i),A)=\sum_{j\in I} \int_{\text{supp}(\nu^{\varepsilon})} \int_{\Theta} \mathbbm{1}_A(w_{\theta}(y)+h,j)\pi_{ij}(w_{\theta}(y)+h)p(y,{\theta})\,d\theta\,\nu^{\varepsilon}(dh)
\end{gather}
for any $(y,i)\in X$, $A\in\mathcal{B}_X$. 
More precisely, \cite[Theorem 4.4]{dawid} says that if $\mu_*\in\mathcal{M}_1(X)$ is an invariant measure of the Markov operator $P$, then $\nu_*:=\mu_* G$ is an invariant  measure of the Markov semigroup $({P}_t)_{t\in\mathbb{R}_+}$, and $\nu_* W = \mu_*$. Conversely, if $\nu_*\in\mathcal{M}_1(X)$ is an invariant measure of $({P}_t)_{t\in\mathbb{R}_+}$, then \hbox{$\mu_*:=\nu_* W$} is an invariant  measure of $P$, and $\mu_* G=\nu_*$.

Finally, let us denote the renewal counting process with arrival times $\tau_n$, $n\in\mathbb{N}_0$, by $(N_t)_{t\in\mathbb{R}_+}$, i.e.
\begin{align}\label{def:N_t}
N_t:=\max\{n\in\n_0: \tau_n\leq t\}\;\;\;\text{for} \;\;\;t\in\mathbb{R}_+.
\end{align}

\section{The main result}

Consider the Markov chain $(X_n)_{n\in\mathbb{N}_0}$ with transition law $P$, given by \eqref{def:Pi_epsilon}, as well as the piecewise deterministic Markov process $(X(t))_{t\in\mathbb{R}_+}$, defined by \eqref{def:X(t)}. Further, recall that under hypotheses \eqref{cnd:a1}-\eqref{cnd:a5} and \eqref{eq:balance1} both the semigroup $(P_t)_{t\in\mathbb{R}_+}$ and the operator $P$ possess unique invariant distributions, denoted by $\nu_*\in\mathcal{M}_1(X)$ and $\mu_*\in\mathcal{M}_1(X)$, respectively. Moreover, we know that $\nu_*=\mu_*G$, where $G$ is defined in \eqref{def_G}.

Let $g\in Lip_b(X)$ be an arbitrary non-constant function, and define $\bar{g}=g-\langle g,\nu_*\rangle$. Following \eqref{def:s_n} and \eqref{def:s_t}, we can introduce 
\begin{gather}\label{def:s_n_X}
s_n(G\bar{g})=\cfrac{\sum_{i=0}^{n-1}G\bar{g}(X_i)}{\sqrt{2n\ln(\ln (n))}}\;\;\;\text{for}\;\;\;n>e,\qquad s_n(G\bar{g})=0\;\;\;\text{for}\;\;\;n\leq e,\\
s(\bar{g})(t)=\cfrac{\int_0^t\bar{g}(X(s))ds}{\sqrt{2t\ln(\ln (t))}} 
\;\;\;\text{for} \;\;\;t>e,\qquad\text{and}\qquad
s(\bar{g})(t)=0\;\;\;\text{for} \;\;\;t\leq e.
\end{gather}

We are now ready to state our main result, whose proof is presented in the ramainder of the paper.

\begin{theorem}\label{MAIN_THM}
Suppose that conditions \eqref{cnd:a1}-\eqref{cnd:a5} hold with the constants satisfying \eqref{eq:balance1}. Then, for any non-constant function $g\in Lip_b(X)$ and any initial measure $\mu\in\mathcal{M}_{1,2+r}^V(X)$ with $V$ given by \eqref{def:V}, the process $(g(X(t)))_{t\in\mathbb{R}_+}$ satisfies the LIL.
\end{theorem}

\subsection{The proof of the main result}
According to the definition introduced in Section \ref{sec:def_lil}, we need to prove that 
\begin{align*}
\limsup\limits_{t\to\infty}s\left(\bar{g}\right)(t)=\overline{\sigma}(\bar{g})
\;\;\;
\text{and}\;\;\;
\liminf\limits_{t\to\infty}s\left(\bar{g}\right)(t)=-\overline{\sigma}(\bar{g})\;\;\;\mathbb{P}\text{-a.s.}
\end{align*}
for some $\overline{\sigma}(\bar{g})\in(0,\infty)$.

Recall that, for any $t\in\mathbb{R}_+$, $N_t$ is given by \eqref{def:N_t}. Further, note that whenever $t \geq \tau_3$, which in other words means that $N_t>e$, we have
\begin{align*}
s(\bar{g})(t)
=\cfrac{\sqrt{2N_t\ln(\ln (N_t))}}{\sqrt{2t\ln(\ln(t))}}\left(\cfrac{1}{\sqrt{2N_t\ln(\ln(N_t))}}\sum_{i=0}^{N_{t}-1}\int_{\tau_i}^{\tau_{i+1}}\bar{g}\left(X(s)\right)ds+R_t(\bar{g})\right),
\end{align*}
where 
\begin{align*}
R_t(\bar{g}):=\cfrac{1}{\sqrt{2N_t\ln(\ln(N_t))}}\int_{\tau_{N_t}}^t\bar{g}\left(X(s)\right)ds.
\end{align*}
We can further write
\begin{align}\label{eq:s_t_3_parts}
\begin{aligned}
s(\bar{g})(t)=&\cfrac{\sqrt{N_t\ln(\ln (N_t))}}{\sqrt{t\ln(\ln(t))}}\Bigg(\cfrac{1}{\sqrt{2N_t\ln(\ln(N_t))}}
\sum_{i=0}^{N_{t}-1}\left(
\int_{\tau_i}^{\tau_{i+1}}\bar{g}\left(X(s)\right)ds-\frac{1}{\lambda}G\bar{g}(X_i)\right)\\
&+R_t(\bar{g})
+\frac{1}{\lambda}s_{N_t}(G\bar{g})\Bigg),
\end{aligned}
\end{align}
where $s_{N_t}(G\bar{g})$ is defined as in \eqref{def:s_n_X}. 
Referring to the elementary renewal theorem, which says that 
\begin{align}\label{renewal_thm}
\lim_{t\to\infty}\frac{N_t}{t}=\lambda\;\;\;\mathbb{P}\text{-a.s.},
\end{align}
we obtain that 
\begin{align}\label{eq:sqrt_lambda}
\lim_{t\to\infty}\cfrac{\sqrt{N_t\ln(\ln(N_t))}}{\sqrt{t\ln(\ln(t))}}=\sqrt{\lambda}\;\;\;\mathbb{P}\text{-a.s.}
\end{align}

For any $t\in\mathbb{R}_+$, let us now introduce the following notation:
\begin{align*}
&I_1(t):=\cfrac{1}{\sqrt{2N_t\ln(\ln(N_t))}}
\sum_{i=0}^{N_{t}-1}\left(
\int_{\tau_i}^{\tau_{i+1}}\bar{g}\left(X(s)\right)ds-\frac{1}{\lambda}G\bar{g}(X_i)\right),\\
&I_2(t):=R_t(\bar{g}),\\
&I_3(t):=\frac{1}{\lambda}s_{N_t}(G\bar{g}).
\end{align*}
The asymptotic behavior of each of these components shall be analyzed separately.

First of all, we have
\begin{align}\label{eq:estim_R_t}
|R_t(\bar{g})|\leq \|\bar{g}\|_{\infty}\cfrac{\Delta\tau_{N_t+1}}{\sqrt{2N_t\ln(\ln(N_t))}}\;\;\;\mathbb{P}\text{-a.s.}\;\;\;\text{for}\;\;\;t\geq\tau_3.
\end{align}
Observe that the right-hand side of the above inequality tends to zero. Indeed, note that 
\[\sum_{n=3}^{\infty}
\mathbb{P}
\left(\cfrac{\Delta\tau_{n+1}}{\sqrt{2n\ln(\ln(n))}}\geq \varepsilon\right)
=\sum_{n=3}^{\infty}e^{-\lambda\varepsilon\sqrt{2n\ln(\ln(n))}}<\infty.\]
Hence, due to the Borel-Cantelli lemma, 
\begin{align*}
\lim_{n\to\infty}\cfrac{\Delta\tau_{n+1}}{\sqrt{2n\ln(\ln(n))}}=0\;\;\;\mathbb{P}\text{-a.s.},
\end{align*}
whence also
\begin{align*}
\lim_{t\to\infty}\cfrac{\Delta\tau_{N_t+1}}{\sqrt{2N_t\ln(\ln(N_t))}}=0\;\;\;\mathbb{P}\text{-a.s.},
\end{align*}
which follows from \eqref{renewal_thm}. 
Finally, referring to \eqref{eq:estim_R_t}, we see that  
\begin{align}\label{eq:I_2}
\lim_{t\to\infty}I_2(t)=0\;\;\;\mathbb{P}\text{-a.s.}
\end{align}

While investigating $I_3$, we shall refer to \cite[Theorem 4.1]{lil_chw}. Note that the Markov chain $(X_n)_{n\in\mathbb{N}_0}$, for which the sequence $(s_{n}(G\bar{g}))_{n\in\mathbb{N}_0}$ is defined, satisfies all the assumptions required in \cite[Theorem 4.1]{lil_chw}. Therefore the only conditions that need to be proven are $G\bar{g}\in Lip_b(X)$ and $\langle G\bar{g},\mu_*\rangle=0$, where the latter follows immediately from the definition of $\bar{g}$ and the fact that $\langle G\bar{g}, \mu_*\rangle=\langle \bar{g}, \nu_*\rangle$ (cf. \cite[Theorem 4.4]{dawid}). 
Since the boundedness of $G\bar{g}$ is also obvious, it remains to show its Lipschitz-continuity. Note that, according to \eqref{cnd:a2}, we have
\begin{align*}
|G\bar{g}(y_1,i_1)&-G\bar{g}(y_2,i_2)|\leq
\int_0^{\infty}\lambda e^{-\lambda t}
\left|g\left(S_{i_1}(t,y_1),i_1\right) - g\left(S_{i_2}(t,y_2),i_2\right)\right|dt\\
&\leq  |g|_{Lip}\int_0^{\infty}\lambda e^{-\lambda t}\left(\left\|S_{i_1}(t,y_1)-S_{i_2}(t,y_2)\right\|+c d(i_1,i_2)\right)dt\\
&\leq  |g|_{Lip}\left(\lambda L\|y_1-y_2\|\int_0^{\infty}e^{-(\lambda-\alpha)t}dt+d(i_1,i_2)\bar{L}\int_0^{\infty}\lambda e^{-\lambda t}t\,dt+c d(i_1,i_2)\right)\\
&=|g|_{Lip}\left(\cfrac{\lambda L}{\lambda-\alpha}\|y_1-y_2\|+d(i_1,i_2)\left(\cfrac{\bar{L}}{\lambda}+c\right)\right)\\
&\leq |g|_{Lip}\left(\cfrac{\lambda L}{\lambda-\alpha}+\cfrac{\bar{L}}{\lambda}+c\right)\varrho_{c}((y_1,i_1),(y_2,i_2))\;\;\;\text{ for any}\;\;\;(y_1,i_1),(y_2,i_2)\in X,
\end{align*}
which guarantees that $G\bar{g}\in Lip_b(X)$. Therefore it follows from \cite[Theorem 4.1]{lil_chw} that  
\begin{align}\label{limit_s_n}
\limsup\limits_{n\to\infty}s_{n}(G\bar{g})=\sigma(G\bar{g})
\;\;\;\text{and}\;\;\;
\liminf\limits_{n\to\infty}s_{n}(G\bar{g})=-\sigma(G\bar{g})\;\;\;\mathbb{P}\text{-a.s.},
\end{align}
where, for any function $h\in Lip_b(X)$,
$$\sigma^2(h)=\mathbb{E}_{\mu_*}\left(\left(\sum_{i=0}^{\infty}P^ih(X_1)-\sum_{i=0}^{\infty}P^ih(X_0)+h(X_0)\right)^2\right),$$
and $\mathbb{E}_{\mu_*}$ is the expected value corresponding to the probability measure $\mathbb{P}_{\mu_*}$ defined on $(\Omega,\mathcal{F})$ such that $\mathbb{P}_{\mu_*}(X_0\in A)=\mu_*(A)$ for $A\in\mathcal{F}$.  
Hence, due to \eqref{limit_s_n} and \eqref{renewal_thm}, we obtain
\begin{gather}\label{eq:I_3}
\begin{aligned}
\limsup\limits_{t\to\infty}I_3(t)=\frac{1}{\lambda}\limsup\limits_{t\to\infty}s_{N_t}(G\bar{g})=\frac{1}{\lambda}\sigma(G\bar{g})\;\;\;&\mathbb{P}\text{-a.s.},\\
\text{and}\;\;\;\liminf\limits_{t\to\infty}I_3(t)=-\frac{1}{\lambda}\sigma(G\bar{g})\;\;\;&\mathbb{P}\text{-a.s.}
\end{aligned}
\end{gather}
Note that $\sigma(G\bar{g})<\infty$, which is explained in details in \cite{lil_chw}. 

Finally, to analyze the asymptotic behaviour of $I_1(t)$, we need to appeal to \hbox{\cite[Theorem 1]{hs}}, whose assertion guarantees the LIL for certain square integrable martingales. Let us first introduce the sequence $(M_n(\bar{g}))_{n\in\mathbb{N}_0}$ given by
\begin{align}\label{def:M_n}
M_0(\bar{g})=0\;\;\;\text{and}\;\;\;M_n(\bar{g})=\sum_{k=0}^{n-1} \left(\int_{\tau_k}^{\tau_{k+1}} \bar{g}\left(X(s)\right)\,ds-\frac{1}{\lambda} G\bar{g}(X_k)\right)\;\;\;\text{for}\;\;\;n\in\n.
\end{align}
Note that $(M_n(\bar{g}))_{n\in\n_0}$ is a martingale with respect to the natural filtration $(\mathcal{F}_n)_{n\in\n_0}$ of $(X_n,\Delta\tau_n)_{n\in\n_0}$. Indeed, we have
\begin{align}\label{def:Z_n}
\begin{aligned}
Z_{n+1}(\bar{g}):=M_{n+1}(\bar{g})-M_n(\bar{g})
&=\int_{\tau_n}^{\tau_{n+1}}\bar{g}\left(X(s)\right)ds-\frac{1}{\lambda}G\bar{g}\left(X_n\right)\\
&=\int_{\tau_n}^{\tau_{n+1}}\bar{g}\left(S_{\xi_n}(s-\tau_n,Y_n),\xi_n\right)ds-\frac{1}{\lambda}G\bar{g}\left(Y_n,\xi_n\right)\\
&=\int_0^{\Delta\tau_{n+1}}\bar{g}\left(S_{\xi_n}(s,Y_n),\xi_n\right)ds-\frac{1}{\lambda}G\bar{g}\left(Y_n,\xi_n\right),
\end{aligned}
\end{align}
whence, appealing to \eqref{properties_P}, for any $(y,i,u)\in X\times \mathbb{R}_+$, we get
\begin{align*}
\mathbb{E}\left(Z_{n+1}(\bar{g})|Y_n=y,\xi_n=i,\Delta\tau_n=u\right)
=&\int_{\mathbb{R}}\int_0^t\bar{g}\left(S_i(s,y),i\right)ds\,\mathbb{P}\left(\Delta\tau_{n+1}\in dt\right)\\
&-\frac{1}{\lambda}G\bar{g}\left(y,i\right)\\
=&\int_0^{\infty}\lambda e^{-\lambda t}\int_0^t\bar{g}\left(S_i(s,y),i\right)ds\, dt-\frac{1}{\lambda}G\bar{g}\left(y,i\right)\\
=&\int_0^{\infty}\int_s^{\infty}\lambda e^{-\lambda t}dt\,\bar{g}\left(S_i(s,y),i\right)ds
-\frac{1}{\lambda}G\bar{g}\left(y,i\right)\\
=&\int_0^{\infty} e^{-\lambda s}\bar{g}\left(S_i(s,y),i\right)ds-\frac{1}{\lambda}G\bar{g}\left(y,i\right)=0,
\end{align*}
which, by the Markov property of the chain $(X_n,\Delta\tau_n)_{n\in\mathbb{N}_0}$, implies that $(M_n(\bar{g}))_{n\in\n_0}$ is a~martingale. 
Further, we also obtain
\begin{align*}
\mathbb{E}\left(Z_{n+1}^2(\bar{g})\right)
\leq &2\mathbb{E}\left(\left(\int_0^{\Delta\tau_{n+1}}\bar{g}\left(S_{\xi_n}(s,Y_n),\xi_n\right)ds\right)^2\right)+2\mathbb{E}\left(\left(\frac{1}{\lambda}G\bar{g}(Y_n,\xi_n)\right)^2\right)\\
\leq & 2\|\bar{g}\|_{\infty}^2\mathbb{E}\left(\left(\Delta\tau_{n+1}\right)^2\right)+\frac{2}{\lambda^2}\|\bar{g}\|_{\infty}^2=\frac{6}{\lambda^2}\|\bar{g}\|_{\infty}^2,
\end{align*}
which means that the martingale increments $Z_n(\bar{g})=M_n(\bar{g})-M_{n-1}(\bar{g})$, $n\in\n$, are uniformly bounded in the $\mathcal{L}^2(\mathbb{P})$-norm, and thus the martingale itself is square-integrable, as required in \cite[Theorem 1]{hs}. 

Now, define
\begin{align*}
h_n^2(\bar{g}):=\mathbb{E}\left(M_n^2(\bar{g})\right)\;\;\;\text{for}\;\;\;n\in\mathbb{N}_0.
\end{align*}
It will be clarified later on (in Section \ref{sec:marting}) that there exists $\bar{n}\in\mathbb{N}$ such that $h_n(\bar{g})>0$ for every $n\geq \bar{n}$. 
We need to establish the following conditions:
\begin{gather}\label{eq:h_n_1}
\lim_{n\to\infty}\frac{1}{h_n^2(\bar{g})}\sum_{l=1}^{n}Z_l^2(\bar{g})=1\;\;\;\mathbb{P}\text{-a.s.},\\
\sum_{n=\bar{n}}^{\infty}
h_n^{-4}(\bar{g})\mathbb{E}\left(Z_n^4(\bar{g})
\mathbbm{1}_{\left\{|Z_n(\bar{g})|<\upsilon h_n(\bar{g})\right\}}\right)
< \infty\;\;\;\text{for every}\;\;\;\upsilon>0,\label{eq:lil1}\\
\sum_{n=\bar{n}}^{\infty}
h_n^{-1}(\bar{g})\mathbb{E}\left(|Z_n(\bar{g})|
\mathbbm{1}_{\left\{|Z_n(\bar{g})|\geq\vartheta h_n(\bar{g})\right\}}\right)
<\infty\;\;\;\text{for every}\;\;\;\vartheta>0,\label{eq:lil2}
\end{gather}
which, in view of \cite[Theorem 1]{hs}, imply the LIL for the martingale $(M_n(\bar{g}))_{n\in\mathbb{N}_0}$. To be more precise, according to \cite[Theorem 1]{hs}, the sequence $(M_n(\bar{g}))_{n\in\mathbb{N}_0}$ satisfies the Strassen invariance principle for the LIL with the normalizing factors $$\frac{1}{\sqrt{2h_n^2(\bar{g})\ln(\ln (h_n^2(\bar{g})))}},\;\;\;n\geq \bar{n}.$$ 
In particular, it also satisfies the LIL itself, which, in this case, means that
\begin{gather*}
\limsup\limits_{n\to\infty}\cfrac{M_n(\bar{g})}{\sqrt{2h_n^2(\bar{g})\ln(\ln (h_n^2(\bar{g})))}}=1\;\;\;\mathbb{P}\text{-a.s.},\;\;\;
\liminf\limits_{n\to\infty}\cfrac{M_n(\bar{g})}{\sqrt{2h_n^2(\bar{g})\ln(\ln (h_n^2(\bar{g})))}}=-1\;\;\;\mathbb{P}\text{-a.s.},
\end{gather*}
and so, according to \eqref{renewal_thm}, we further obtain
\begin{gather*}
\limsup\limits_{t\to\infty}\cfrac{M_{N_t}(\bar{g})}{\sqrt{2h_{N_t}^2(\bar{g})\ln(\ln (h_{N_t}^2(\bar{g})))}}=1\;\;\;\mathbb{P}\text{-a.s.},\;\;\;
\liminf\limits_{n\to\infty}\cfrac{M_{N_t}(\bar{g})}{\sqrt{2h_{N_t}^2(\bar{g})\ln(\ln (h_{N_t}^2(\bar{g})))}}=-1\;\;\;\mathbb{P}\text{-a.s.}
\end{gather*}
Let the part of the proof in which we verify \eqref{eq:h_n_1}-\eqref{eq:lil2} be postponed into the subsequent section, namely Section \ref{sec:marting}, in which we shall also prove that 
\begin{align}\label{eq:conv_to_sigma}
\lim_{t\to\infty}\cfrac{\sqrt{h_{N_t}^2(\bar{g})\ln\left(\ln(h_{N_t}^2(\bar{g}))\right)}}{\sqrt{N_t\ln(\ln(N_t))}}=\tilde{\sigma}(\bar{g})\;\;\;\mathbb{P}\text{-a.s.},
\end{align}
where
\begin{align}\label{def:tilde_sigma}
\tilde{\sigma}^2(\bar{g}):=\mathbb{E}_{\mu_*}\left(Z_1^2(\bar{g})\right)=\mathbb{E}_{\mu_*}\left(M_1^2(\bar{g})\right)\in(0,\infty).
\end{align}
Then, provided that \eqref{eq:h_n_1}-\eqref{eq:lil2} and \eqref{eq:conv_to_sigma} are established, we obtain
\begin{align}\label{eq:I_1}
\begin{aligned}
\limsup\limits_{t\to\infty}I_1(t)=\tilde{\sigma}(\bar{g})\;\;\;\text{and}\;\;\;
\liminf\limits_{t\to\infty}I_1(t)=-\tilde{\sigma}(\bar{g})\;\;\;\mathbb{P}\text{-a.s.}
\end{aligned}
\end{align}

Finally, combining \eqref{eq:s_t_3_parts} with \eqref{eq:sqrt_lambda}, \eqref{eq:I_2}, \eqref{eq:I_3} and \eqref{eq:I_1}, we obtain
\begin{align}
\limsup\limits_{t\to\infty}s(\bar{g})(t)=\overline{\sigma}(\bar{g})
\;\;\;\text{and}\;\;\;
\liminf\limits_{t\to\infty}s(\bar{g})(t)=-\overline{\sigma}(\bar{g})
\;\;\;\mathbb{P}\text{-a.s.},
\end{align}
where 
\begin{align*}
\overline{\sigma}(\bar{g}):=\sqrt{\lambda}\left(\frac{1}{\lambda}\sigma(G\bar{g})+\tilde{\sigma}(\bar{g})\right)\in(0,\infty).
\end{align*}
The proof of Theorem \ref{MAIN_THM} is therefore completed (provided that \eqref{eq:h_n_1}-\eqref{eq:conv_to_sigma} are established, which shall be done in the upcoming section).

\subsection{The proof of the LIL for the appropriate martingale}\label{sec:marting}

Let us consider
\[\mathcal{Z}:=\{((x_1,t),(x_2,s))\in\left(X\times\mathbb{R}_+\right)^2:\,t=s\},\]
and, for any $A\in\mathcal{B}_{X^2}$, define
\begin{align*}
(A)_{\mathcal{Z}}:=\left\{\left(\left(x_1,t\right),\left(x_2,t\right)\right)\in\mathcal{Z}:\,\left(x_1,x_2\right)\in A\right\}.
\end{align*}
Further, introduce
$\widetilde{Q}:\mathcal{Z}\times\mathcal{B}_{\mathcal{Z}}\to[0,1]$ given by
\begin{align}\label{def:tilde_Q}
\begin{aligned}
\widetilde{Q}\left(\left(\left(x_1,s\right),\left(x_2,s\right)\right),B\right)
=&\int_{\text{supp}(\nu^{\varepsilon})}
\int_0^{\infty}\lambda e^{-\lambda t}\int_{\Theta}
\Bigg(\sum_{j\in I}\mathbbm{1}_B\left({\mathbf{w}}_j\left(x_1,x_2,t,\theta,h\right)\right)\\
&\times\boldsymbol{\pi}_j\left(x_1,x_2,t,\theta,h\right)\Bigg)\mathbf{p}\left(x_1,x_2,t,\theta\right)\,d\theta\,dt\,\nu^{\varepsilon}(dh)
\end{aligned}
\end{align}
for $((x_1,s),(x_2,s))\in \mathcal{Z}$ and $B\in\mathcal{B}_{\mathcal{Z}}$ such that $x_1=(y_1,i_1)$, $x_2=(y_2,i_2)$, where 
\begin{gather*}
\mathbf{w}_j\left(x_1,x_2,t,\theta,h\right)
=\left(\left(w_{\theta}\left(S_{i_1}(t,y_1)+h\right),j,t\right),\left(w_{\theta}\left(S_{i_2}(t,y_2)+h\right),j,t\right)\right),\\
\boldsymbol{\pi}_j\left(x_1,x_2,t,\theta,h\right)=\pi_{i_1,j}\left(w_{\theta}\left(S_{i_1}(t,y_1)\right)+h\right)\wedge \pi_{i_2,j}\left(w_{\theta}\left(S_{i_2}(t,y_2)\right)+h\right),\\
\mathbf{p}\left(x_1,x_2,t,\theta\right)=p(S_{i_1}(t,y_1),\theta)\wedge p(S_{i_2}(t,y_2),\theta).
\end{gather*}
Note that $\widetilde{Q}$ is a substochastic kernel, and, for any $x_1,x_2\in X$, $t\in\mathbb{R}_+$, $B\in\mathcal{B}_X$, satisfies the following properties:
\begin{gather*}
\widetilde{Q}\left(\left(\left(x_1,t\right),\left(x_2,t\right)\right),\left(B\times X\right)_{\mathcal{Z}}\right)\leq \Pi\left(\left(x_1,t\right),B\times\mathbb{R}_+\right),\\
\widetilde{Q}\left(\left(\left(x_1,t\right),\left(x_2,t\right)\right),\left(X\times B\right)_{\mathcal{Z}}\right)\leq \Pi\left(\left(x_2,t\right),B\times\mathbb{R}_+\right).
\end{gather*}
For any given distribution $\textbf{m} \in \mathcal{M}_1(X^2)$, on the coordinate space $(\widetilde{\Omega},\widetilde{F})$  associated with $\mathcal{Z}$, we can now construct a probability measure $\widetilde{\mathbb{C}}$ so that 
$$\widetilde{\mathbb{C}}\left( \left(\widetilde{X}_0^{(1)}, \widetilde{X}_0^{(2)}\right) \in A, \, \widetilde{\Delta\tau}_0=0\right) = \textbf{m}(A) \;\;\; \text{for any} \;\;\;A\in \mathcal{B}_{X^2},$$
and the canonical Markovian coupling $((\widetilde{X}_n^{(1)},{\widetilde{\Delta\tau}}_n),(\widetilde{X}_n^{(2)},{\widetilde{\Delta\tau}}_n))_{n\in\mathbb{N}_0}$ of $\Pi$, defined on this space, is governed by the transition probability kernel of the form 
$$\widetilde{C}=\widetilde{Q}+\widetilde{R},$$
where $\widetilde{Q}$ is defined by \eqref{def:tilde_Q}, and $\widetilde{R}$ stands for a complementary substochastic kernel on $\mathcal{Z}\times B_{\mathcal{Z}}$. The latter can be specified by defining the corresponding family of measures on rectangles $\{A\times B: A,B\in \mathcal{B}_X\}$ as follows:
\begin{align*}
\widetilde{R}\left(\left(\left(x_1,t\right),\left(x_2,t\right)\right),\left(A\times B\right)_{\mathcal{Z}}\right)
=&\frac{1}{1-\widetilde{Q}\left(\left(\left(x_1,t\right),\left(x_2,t\right)\right),\mathcal{Z}\right)}\\
&\times
\left(\Pi\left(\left(x_1,t\right),A\right)-\widetilde{Q}\left(\left(\left(x_1,t\right),\left(x_2,t\right)\right),\left(A\times X\right)_{\mathcal{Z}}\right)\right)\\
&\times\left(\Pi\left(\left(x_2,t\right),B\right)-\widetilde{Q}\left(\left(\left(x_1,t\right),\left(x_2,t\right)\right),\left(X\times B\right)_{\mathcal{Z}}\right)\right),
\end{align*}
when $\widetilde{Q}(((x_1,t),(x_2,t)),\mathcal{Z})<1$, and $\widetilde{R}(((x_1,t),(x_2,t)),(A\times B)_{\mathcal{Z}})=0$ otherwise.

Now, define $Q:X^2\times\mathcal{B}_{X^2}\to[0,1]$ and $C:X^2\times\mathcal{B}_{X^2}\to[0,1]$ as the kernels  which, for any $(x_1,x_2)\in X^2$, $t\in\mathbb{R}_+$ and $A\in\mathcal{B}_{X^2}$, satisfy
\begin{align}\nonumber
Q\left(\left(x_1,x_2\right),A\right)=\widetilde{Q}\left(\left(\left(x_1,0\right),\left(x_2,0\right)\right),\left(A\right)_{\mathcal{Z}}\right)=\widetilde{Q}\left(\left(\left(x_1,t\right),\left(x_2,t\right)\right),\left(A\right)_{\mathcal{Z}}\right),\\
\label{def:couplingC} 
C\left(\left(x_1,x_2\right),A\right)=\widetilde{C}\left(\left(\left(x_1,0\right),\left(x_2,0\right)\right),\left(A\right)_{\mathcal{Z}}\right)=\widetilde{C}\left(\left(\left(x_1,t\right),\left(x_2,t\right)\right),\left(A\right)_{\mathcal{Z}}\right).
\end{align}

Later on in this paper, we will write $\widetilde{\mathbb{E}}_{x_1,x_2}$ for the expected value corresponding to the measure 
\[\widetilde{\mathbb{C}}_{x_1,x_2}:=\widetilde{\mathbb{C}}\left(\cdot\Big|\widetilde{X}^{(1)}_0=x_1,\widetilde{X}^{(2)}_0=x_2\right),\;\;\;x_1,x_2\in X.\]

Let us indicate that the model under consideration enjoys all the hypotheses assumed in \hbox{\cite[Theorem 2.1]{ks}} (see the proof of \cite[Theorem 4.1]{dawid}, where these conditions are verified), which, in particular, means that
\begin{itemize}
\item[(B0)] \phantomsection \label{cnd:B0} The Markov operator $P$ is Feller.
\item[(B1)] \phantomsection \label{cnd:B1} There exist constants $a\in (0,1)$ and $b\in(0,\infty)$ such that
$$PV(x)\leq aV(x)+b\;\;\;\text{for every}\;\;\;x\in X,$$
where $V$ is given by \eqref{def:V}.
\end{itemize}
Moreover, letting 
$$F=\left\{\left(\left(y_1,i_1\right),\left(y_2,i_2\right)\right)\in X^2:\;i_1=i_2\right\}\cup\left\{\left(x_1,x_2\right)\in X^2:\;V(x_1)+V(x_2)<\frac{4b}{1-a}\right\},$$
the following statements hold:
\begin{itemize}
\item[(B2)] \phantomsection \label{cnd:B2} 
We have 
$\,\text{supp}\,Q(x,y,\cdot)\subset F\,$ 
and
\begin{align*}
\int_{X^2}\varrho(u,v)\,Q(x,y,du\times dv)\leq \beta \varrho(x,y)\;\;\;\text{for any}\;\;\;(x,y)\in F\;\;\;\text{and some}\;\;\;\beta\in(0,1).
\end{align*}
\item[(B3)] \phantomsection \label{cnd:B3} Letting $U(r):=\{(u,v)\in F:\varrho(u,v)\leq r\}$ for any  $r>0$,  we have
$$\inf_{(x,y)\in F}Q(x,y,U\left(\beta\varrho(x,y)\right))>0.$$
\item[(B4)] \phantomsection \label{cnd:B4}  There exists $l>0$ such that
$$Q\left(x,y,X^2\right)\geq 1-l\varrho(x,y)\quad\text{for every}\quad (x,y)\in F.$$
\item[(B5)] \phantomsection \label{cnd:B5} 
There exist $\gamma\in(0,1)$ and $\hat{c}>0$ such that
$$\widetilde{\mathbb{E}}_{x_1,x_2}(\gamma^{-\rho})\leq \hat{c},\quad\text{whenever}\quad V(x)+V(y)<4b(1-a)^{-1},$$
where $V$ is given by \eqref{def:V} and
\begin{align}\label{def:rho}
\rho
=\inf\left\{n\in\mathbb{N}:\;\left(\widetilde{X}^{(1)}_n,\widetilde{X}^{(2)}_n\right)\in F\;\text{and}\;V\left(\widetilde{X}^{(1)}_n\right)+V\left(\widetilde{X}^{(2)}_n\right)<\frac{4b}{1-a}\right\}.
\end{align}
\end{itemize}

For 
$(Z_n(\bar{g}))_{n\in\mathbb{N}}$, given by 
\eqref{def:Z_n}, 
let us consider the sequences of their copies 
$(\widetilde{Z}_n^{(i)}(\bar{g}))_{n\in\mathbb{N}}$, $i\in\{1,2\}$, defined on $(\widetilde{\Omega},\widetilde{\mathcal{F}},\widetilde{\mathbb{C}})$ as follows:
\begin{gather}\label{def:tilde_Z_n}
\widetilde{Z}_n^{(i)}(\bar{g})=Z_n(\bar{g})\left(\left(\widetilde{X}_0^{(i)},\widetilde{\Delta\tau}_0\right), \left(\widetilde{X}_1^{(i)},\widetilde{\Delta\tau}_1\right), \ldots\right)\;\;\;\text{for}\;\;\;n\in\mathbb{N}_0\;\;\;\text{and}\;\;\; i\in\{1,2\}.
\end{gather}

According to \cite[Lemmas 3.4 and 3.5]{lil_chw}, we can now state the following result.
\begin{lemma}
Suppose that  
\begin{align}\label{cond:a}
\sum_{n=1}^{\infty}\widetilde{\mathbb{E}}_{x_1,x_2}|\widetilde{Z}_n^{(1)}(\bar{g})-\widetilde{Z}_n^{(2)}(\bar{g})|<\infty\;\;\;\text{for all}\;\;\;x_1,x_2\in X,
\end{align}
and that there exists $r\in(0,2)$ such that
\begin{align}\label{cond:b}
\sup_{n\in\mathbb{N}}{\mathbb{E}}
|{Z}_n(\bar{g})|^{2+r}<\infty\;\;\;\text{for any}\;\;\;i\in\{1,2\}.
\end{align}
Then 
\[\lim_{n\to\infty}\cfrac{h_n^2(\bar{g})}{n}=\tilde{\sigma}^2(\bar{g})\]
with $\tilde{\sigma}(\bar{g})$ given by \eqref{def:tilde_sigma}, which further yields
\[\lim_{n\to\infty}\cfrac{\sqrt{h_n^2(\bar{g})\ln\left(\ln\left(h_n^2(\bar{g})\right)\right)}}{\sqrt{n\ln(\ln(n))}}=\tilde{\sigma}(\bar{g}),\]
and consequently \eqref{eq:conv_to_sigma} holds.
Moreover, conditions \eqref{cond:a}, \eqref{cond:b} imply that there exists $\bar{n}\in\n$ such that $h_n(\bar{g})>0$ for all $n\geq \bar{n}$, and  that hypotheses 
\eqref{eq:h_n_1}-\eqref{eq:lil2} hold. Hence, due to \cite[Theorem 1]{hs}, the martingale $(M_n(\bar{g}))_{n\in\mathbb{N}_0}$, given by \eqref{def:M_n}, satisfies the LIL. 
\end{lemma}

In view of the above lemma, to finalise the proof of Theorem 3.1, it remains to establish hypotheses \eqref{cond:a}-\eqref{cond:b}.

Let us introduce the function $F(\bar{g}):X\times \mathbb{R}_+\to\mathbb{R}_+$ given by
\begin{align}\label{def:F(g)}
F(\bar{g})(x,t)=\int_0^t\bar{g}\left(S_i(s,y),i\right)ds \;\;\;\text{for any}\;\;\;x=(y,i)\in X,\;t\in\mathbb{R}_+.
\end{align}
We then have
\begin{align}\label{eq:E_x1,x2}
\begin{aligned}
\widetilde{\mathbb{E}}_{x_1,x_2}
\left|\widetilde{Z}_{n+1}^{(1)}(\bar{g})-\widetilde{Z}_{n+1}^{(2)}(\bar{g})\right|
\leq & \widetilde{\mathbb{E}}_{x_1,x_2}
\left|F(\bar{g})\left(\widetilde{X}_n^{(1)},\widetilde{\Delta\tau}_{n+1}\right)-F(\bar{g})\left(\widetilde{X}_n^{(2)},\widetilde{\Delta\tau}_{n+1}\right)\right|\\
&+\frac{1}{\lambda}\widetilde{\mathbb{E}}_{x_1,x_2}\left|G\bar{g}\left(\widetilde{X}_n^{(1)}\right)-G\bar{g}\left(\widetilde{X}_n^{(2)}\right)\right|.
\end{aligned}
\end{align}
Let us estimate each component on the right hand side of \eqref{eq:E_x1,x2} separately. First of all, according to \eqref{def:couplingC} and \eqref{properties_P}, we have
\begin{align*}
&\widetilde{\mathbb{E}}_{x_1,x_2}
\left|F(\bar{g})\left(\widetilde{X}_n^{(1)},\widetilde{\Delta\tau}_{n+1}\right)-F(\bar{g})\left(\widetilde{X}_n^{(2)},\widetilde{\Delta\tau}_{n+1}\right)\right|\\
&=\int_{X^2}\left(\int_0^{\infty}\lambda e^{-\lambda  t}|F(\bar{g})(u,i,t)-F(\bar{g})(v,j,t)|\,dt\right){C}^n\left(\left(x_1,x_2\right),\left(du\times di\right)\times \left(dv\times dj\right)\right).
\end{align*}
Further, according to \eqref{def:F(g)}, we get
\begin{align*}
\int_0^{\infty}\lambda e^{-\lambda t}|F(\bar{g})(u,i,t)&-F(\bar{g})(v,j,t)|\,dt\\
&\leq \int_0^{\infty}\lambda e^{-\lambda t}\int_0^t\left|\bar{g}\left(S_i(s,u),i\right)-\bar{g}\left(S_j(s,v),j\right)\right|\,ds \,dt\\
&= \int_0^{\infty}\int_s^{\infty}\lambda e^{-\lambda t}\left|\bar{g}\left(S_i(s,u),i\right)-\bar{g}\left(S_j(s,v),j\right)\right|\,dt\,ds\\
&=\int_0^{\infty}\left(\int_s^{\infty}\lambda e^{-\lambda t}dt\right)\left|\bar{g}\left(S_i(s,u),i\right)-\bar{g}\left(S_j(s,v),j\right)\right|\,ds\\
&=\int_0^{\infty}e^{-\lambda s}
\left|\bar{g}\left(S_i(s,u),i\right)-\bar{g}\left(S_j(s,v),j\right)\right|\,ds,
\end{align*}
and therefore
\begin{align}\label{eq:estim_Fg}
\begin{aligned}
&\widetilde{\mathbb{E}}_{x_1,x_2}
\left|F(\bar{g})\left(\widetilde{X}_n^{(1)},\widetilde{\Delta\tau}_{n+1}\right)-F(\bar{g})\left(\widetilde{X}_n^{(2)},\widetilde{\Delta\tau}_{n+1}\right)\right|\\
&\leq\int_{X^2}\int_0^{\infty}e^{-\lambda s}
\left|\bar{g}\left(S_i(s,u),i\right)-\bar{g}\left(S_j(s,v),j\right)\right|\,ds\,{C}^n\left(\left(x_1,x_2\right),\left(du\times di\right)\times \left(dv\times dj\right)\right).
\end{aligned}
\end{align}
The second component on the right-hand side of \eqref{eq:E_x1,x2} can be estimated similarly, i.e.
\begin{align}\label{eq:estim_Gg}
\begin{aligned}
&\frac{1}{\lambda}\widetilde{\mathbb{E}}_{x_1,x_2}\left|G\bar{g}\left(\widetilde{X}_n^{(1)}\right)-G\bar{g}\left(\widetilde{X}_n^{(2)}\right)\right|\\
&\leq\int_{X^2}\int_0^{\infty} e^{-\lambda s}\left|\bar{g}\left(S_i(s,u),i\right)-\bar{g}\left(S_j(s,v),j\right)\right|\,ds\,{C}^n\left(\left(x_1,x_2\right),\left( du\times di\right)\times\left( dv\times dj\right)\right).
\end{aligned}
\end{align}
Combining \eqref{eq:E_x1,x2}, \eqref{eq:estim_Fg} and \eqref{eq:estim_Gg}, we obtain
\begin{align}\label{estim_E_x1_x2}
\begin{aligned}
&\widetilde{\mathbb{E}}_{x_1,x_2}
\left|\widetilde{Z}_n^{(1)}(\bar{g})-\widetilde{Z}_n^{(2)}(\bar{g})\right|\\
&\leq 2\int_{X^2}\int_0^{\infty} e^{-\lambda s}\left|\bar{g}\left(S_i(s,u),i\right)-\bar{g}\left(S_j(s,v),j\right)\right|\,ds\,{C}^n\left(\left(x_1,x_2\right),\left(du\times di\right)\times\left( dv\times dj\right)\right).
\end{aligned}
\end{align}

Consider $\widehat{\mathcal{Z}}=\widehat{\mathcal{Z}}_Q\cup\widehat{\mathcal{Z}}_R$, 
where $\widehat{\mathcal{Z}}_Q:=\mathcal{Z}\times\{1\}$ and $\widehat{\mathcal{Z}}_R:=\mathcal{Z}\times \{0\}$.  
There exists then some probability space $(\widehat{\Omega},\widehat{\mathcal{F}},\widehat{\mathbb{C}})$, on which we can construct a~time-homogeneous canonical Markov chain $((\widehat{X}_n^{(1)},\widehat{\Delta\tau}_n),(\widehat{X}_n^{(2)},\widehat{\Delta\tau}_n),\zeta_n)_{n\in\mathbb{N}_0}$ with $\widehat{\Delta\tau}_0=0$ and $\zeta_0=0$, evolving on $\widehat{\mathcal{Z}}$, and such that its transition probability function $\widehat{C}$ is given by
\begin{align*}
\widehat{C}\big(\left(\left(x_1,t\right),\left(x_2,t\right),\zeta\right),A\big)
=&
\left(\widetilde{Q}
\left(
\left(\left(x_1,t\right),\left(x_2,t\right)\right),
\cdot\right)
\otimes\delta_1\right)(A)
\\
&+\left(\widetilde{R}\left(\left(\left(x_1,t\right),\left(x_2,t\right)\right),\cdot\right)\otimes\delta_0\right)(A)
\end{align*}
for $((x_1,t),(x_2,t),\zeta)\in\widehat{\mathcal{Z}}$ and $A\in\mathcal{B}_{\widehat{\mathcal{Z}}}$ (cf. e.g. \cite{dawid,clt_chw,ks}).  
By convention, we will further write $\widehat{\mathbb{C}}_{x_1,x_2}(\cdot)$ for $\widehat{\mathbb{C}}(\cdot|\widehat{X}_0^{(1)}=x_1,\widehat{X}_0^{(2)}=x_2)$, and we will denote the corresponding expected value by $\widehat{\mathbb{E}}_{x_1,x_2}$, $x_1,x_2\in X$.

Let $\rho$ be given by \eqref{def:rho}, and, for $N\in\mathbb{N}$, define
\begin{align*}
\rho_N:=\inf\left\{n\geq N:\;\left(\widehat{X}_n^{(1)},\widehat{X}_n^{(2)}\right)\in F\;\;\;\text{and}\;\;\;V\left(\widehat{X}_n^{(1)}\right)+V\left(\widehat{X}_n^{(2)}\right)<\frac{4b}{1-a}\right\}.
\end{align*}
Moreover, introduce 
\begin{align*}
\tau:=\inf\left\{n\in\mathbb{N}:\;\left(\left(\widehat{X}_k^{(1)},\widehat{\Delta\tau}_k\right),\left(\widehat{X}_k^{(2)},\widehat{\Delta\tau}_k\right),\zeta_k\right)\in\widehat{\mathcal{Z}}_Q\;\;\;\text{for all}\;\;\;k\geq n\right\},
\end{align*}
and
\begin{align*}
H_{N,n}=\bigcap_{j=N}^n\left\{\zeta_j=1\right\}\;\;\;\text{for}\;\;\;n,N\in\mathbb{N}\;\;\;\text{such that}\;\;\; n>N.
\end{align*}
Note that
\begin{align}\label{eq:note}
\widehat{\mathbb{C}}_{x_1,x_2}
\left(\widehat{\Omega}\backslash H_{N,n}\right)
=\widehat{\mathbb{C}}_{x_1,x_2}
\left(\bigcup_{j=N}^n\{\zeta_j=0\}\right)
\leq \widehat{\mathbb{C}}_{x_1,x_2}(\tau>N)\;\;\;\text{for}\;\;\;n>N,\;\;\;n,N\in\mathbb{N}.
\end{align}

Now, fix $n,N,M$ such that $n>M>N$ and introduce
\begin{align*}
\widehat{\mathbb{C}}_{x_1,x_2}^{n,M,N}(\cdot):=\widehat{\mathbb{C}}_{x_1,x_2}\left(\cdot\cap\left\{\rho_N\leq M\right\}\cap H_{N,n} \right).
\end{align*}
Following the reasoning presented e.g. in \cite{clt_chw}, and applying the estimate \eqref{eq:note}, we obtain
\begin{align*}
\widehat{\mathbb{C}}_{x_1,x_2}(\cdot)\leq \widehat{\mathbb{C}}_{x_1,x_2}^{n,M,N}(\cdot)+\widehat{\mathbb{C}}_{x_1,x_2}\left(\cdot\cap \left\{\rho_N> M\right\}\right)+\widehat{\mathbb{C}}_{x_1,x_2}\left(\cdot\cap \left\{\tau>N\right\}\right),
\end{align*}
and therefore, using \eqref{estim_E_x1_x2} and referring to the fact that $\bar{g}\in Lip_{b}(X)$, we get
\begin{align}\label{estim1}
\begin{aligned}
\widetilde{\mathbb{E}}_{x_1,x_2}
\left|\widetilde{Z}_n^{(1)}(\bar{g})-\widetilde{Z}_n^{(2)}(\bar{g})\right|
\leq & 2|\bar{g}|_{Lip}\int_{X^2}\left(\int_0^{\infty}e^{-\lambda s}\varrho_{c}\left(\left(S_i(s,u),i\right),\left(S_j(s,v),j\right)\right)ds\right)\\
&\times\widehat{\mathbb{C}}_{x_1,x_2}^{n,M,N}\left(\widehat{X}_n^{(1)}\in du\times di,\widehat{X}_n^{(2)}\in dv\times dj\right)\\
&+\frac{4\|\bar{g}\|_{\infty}}{\lambda}\left({\widehat{\mathbb{C}}}_{x_1,x_2}\left( \rho_N> M\right)+\widehat{\mathbb{C}}_{x_1,x_2}\left( \tau>N\right)\right),
\end{aligned}
\end{align}
where $\varrho_{c}$ is given by \eqref{def:rho_c}. 
Further, condition \eqref{cnd:a2}  implies the following:
\begin{align}\label{estim2}
\begin{aligned}
&\int_0^{\infty}e^{-\lambda s}\varrho_{c}\left(\left(S_i(s,u),i\right),\left(S_j(s,v),j\right)\right)ds
\leq \int_0^{\infty}e^{-\lambda s} \left(\left\|S_i(s,u)-S_j(s,v)\right\|+c d(i,j)\right)ds\\
&\hspace{4cm}\leq \int_0^{\infty}e^{-\lambda s}\left(Le^{\alpha s}\|u-v\|+s\bar{L}d(i,j)+c  d(i,j)\right)ds\\
&\hspace{4cm}=L\|u-v\|\int_0^{\infty}e^{-(\lambda-\alpha)s}ds+d(i,j)\int_0^{\infty}\left(\bar{L}se^{-\lambda s}+c e^{-\lambda s}\right)ds\\
&\hspace{4cm}=\frac{L}{\lambda-\alpha}\|u-v\|+\left(\frac{\bar{L}}{\lambda^2}+\frac{c}{\lambda}\right)d(i,j)\\
&\hspace{4cm}\leq\left(\frac{L}{\lambda-\alpha}+\frac{\bar{L}}{\lambda^2}+\frac{1}{\lambda}\right)\varrho_{c}\left((u,i),(v,j)\right).
\end{aligned}
\end{align}
Note that the last inequality holds, since $c$ is required to be sufficiently large. According to \eqref{estim1} and \eqref{estim2}, we obtain
\begin{align}\label{estim3}
\begin{aligned}
\widetilde{\mathbb{E}}_{x_1,x_2}
\left|\widetilde{Z}_n^{(1)}(\bar{g})-\widetilde{Z}_n^{(2)}(\bar{g})\right|
\leq & 2|\bar{g}|_{Lip}\left(\frac{L}{\lambda-\alpha}+\frac{\bar{L}}{\lambda^2}+\frac{1}{\lambda}\right)
\int_{X^2}
\varrho_{c}\left((u,i),(v,j)\right)\\
&\times
\widehat{\mathbb{C}}_{x_1,x_2}^{n,M,N}\left(\widehat{X}_n^{(1)}\in du\times di,\widehat{X}_n^{(2)}\in dv\times dj\right)\\
&+\frac{4\|\bar{g}\|_{\infty}}{\lambda}\left(\widehat{\mathbb{C}}_{x_1,x_2}\left( \rho_N> M\right)+\widehat{\mathbb{C}}_{x_1,x_2}\left( \tau>N\right)\right).
\end{aligned}
\end{align}
Due to \cite[Lemma 2.2]{clt_chw}, there exist constants $c_1,c_2,c_3> 0$, $q_1,q_2,q_3\in(0,1)$ 
 and $p\geq 1$ such that, for any $x_1,x_2\in X$ and $n,N,M\in\mathbb{N}$ satisfying $n>N>M$, the following inequalities hold: 
\begin{gather*}
\int_{X^2}\varrho_{c}((u,i),(v,j))\,\widehat{\mathbb{C}}_{x_1,x_2}^{n,M,N}\left(\widehat{X}_n^{(1)}\in du\times di,\widehat{X}_n^{(2)}\in dv\times dj\right)\leq c_1q_1^{n-M},\\
\widehat{\mathbb{C}}_{x_1,x_2}\left( \rho_N> M\right)\leq c_2q_2^{M-pN}\left(1+V(x_1)+V(x_2)\right),\\
\widehat{\mathbb{C}}_{x_1,x_2}\left( \tau>N\right)\leq c_3q_3^N\left(1+V(x_1)+V(x_2)\right),
\end{gather*}
which, together with \eqref{estim3}, imply
\begin{align*}
\begin{aligned}
\widetilde{\mathbb{E}}_{x_1,x_2}
\left|\widetilde{Z}_n^{(1)}(\bar{g})-\widetilde{Z}_n^{(2)}(\bar{g})\right|
\leq &2|\bar{g}|_{Lip}\left(\frac{L}{\lambda-\alpha}+\frac{\bar{L}}{\lambda^2}+\frac{1}{\lambda}\right)c_1q_1^{n-M}\\
&+\frac{4\|\bar{g}\|_{\infty}}{\lambda}\left(c_2q_2^{M-pN}+c_3q_3^N\right)\left(1+V(x_1)+V(x_2)\right)\\
&\leq C\|\bar{g}\|_{\infty}\left(q_1^{n-M}+q_2^{M-pN}+q_3^N\right)\left(1+V(x_1)+V(x_2)\right)
\end{aligned}
\end{align*}
with 
\[C:=2c_1\left(\frac{L}{\lambda-\alpha}+\frac{\bar{L}}{\lambda^2}+\frac{1}{\lambda}\right)+\frac{4}{\lambda}(c_2+c_3).\]
Now, define $n_0=\lceil 4p\rceil$ and fix an~arbitrary $n>n_0$. Letting $N=\lfloor n/(4p)\rfloor$ and $M=\lceil n/2\rceil$, we obtain
\begin{align*}
\widetilde{\mathbb{E}}_{x_1,x_2}
\left|\widetilde{Z}_n^{(1)}(\bar{g})-\widetilde{Z}_n^{(2)}(\bar{g})\right|
\leq \bar{C}\|\bar{g}\|_{BL}q^n\left(1+V(x_1)+V(x_2)\right)\;\;\;\text{for every}\;\;\;x_1,x_2\in X,
\end{align*}
where $\bar{C}:=C\max\{q_1^{-1},q_2^{-p}\}$ and $q:=\max\{q_1^{1/2},q_2^{1/4},q_3^{1/(4p)}\}\in(0,1)$. Since $\bar{g}$ is bounded, the above estimation also holds (with some $\hat{C}$ in the place of $\bar{C}$) for $n\leq n_0$. We finally get
\begin{align*}
\sum_{n=1}^{\infty}\widetilde{\mathbb{E}}_{x_1,x_2}
\left|\widetilde{Z}_n^{(1)}(\bar{g})-\widetilde{Z}_n^{(2)}(\bar{g})\right|
<\infty\;\;\;\text{for every}\;\;\;x_1,x_2\in X,
\end{align*}
which proves \eqref{cond:a}. 

It now remains to establish \eqref{cond:b}. Referring to \eqref{def:tilde_Z_n}, \eqref{def:Z_n} and \eqref{def_G}, for every $n\in\mathbb{N}$ and any $i\in\{1,2\}$, we obtain
\begin{align*}
{\mathbb{E}}
\left|{Z}_n(\bar{g})\right|^{2+r}
&={\mathbb{E}}\left|
\int_0^{{\Delta\tau}_{n+1}}\bar{g}\left(S_{{\xi}_n}\left(s,{Y}_n\right),{\xi}_n\right)ds-\frac{1}{\lambda}G\bar{g}\left({X}_n\right)\right|^{2+r}\\
=&{\mathbb{E}}\left|
\int_0^{{\Delta\tau}_{n+1}}\bar{g}\left(S_{{\xi}_n}\left(s,{Y}_n\right),{\xi}_n\right)ds-\int_0^{\infty}e^{-\lambda s}\bar{g}\left(S_{{\xi}_n}\left(s,{Y}_n\right),{\xi}_n\right)ds\right|^{2+r}.
\end{align*}
Since $\bar{g}$ is bounded, we further get
\begin{align*}
{\mathbb{E}}
\left|{Z}_n(\bar{g})\right|^{2+r}
\leq \|\bar{g}\|_{\infty}^{2+r}{\mathbb{E}}\left({\Delta\tau}_{n+1}+\frac{1}{\lambda}\right)^{2+r}\;\;\;\text{for every}\;\;\;n\in\n,\;i\in\{1,2\}.
\end{align*}
One can easily prove that, for $r>0$, there exists some $\kappa\in(2,\infty)$ such that 
\begin{align*}
(\psi_1+\psi_2)^{2+r}\leq \kappa\left(\psi_1^{2+r}+\psi_2^{2+r}\right)\;\;\;\text{for any}\;\;\;\psi_1,\psi_2\geq 0,
\end{align*}
whence
\begin{align*}
{\mathbb{E}}
\left|{Z}_n(\bar{g})\right|^{2+r}
\leq \kappa\|\bar{g}\|_{\infty}^{2+r}\left({\mathbb{E}}\left(\left({\Delta\tau}_{1}\right)^{2+r}\right)+\frac{1}{\lambda^{2+r}}\right)\;\;\;\text{for every}\;\;\;n\in\n,\;i\in\{1,2\},
\end{align*}
which is finite, due to the fact that $\widetilde{\Delta\tau}_{1}$ has the exponential distribution. Finally, we get 
\[\sup_{n\in\n}{\mathbb{E}}\left|{Z}_n(\bar{g})\right|^{2+r}<\infty\;\;\;\text{for any}\;\;\;i\in\{1,2\},\]
and the proof is completed.

\section{Acknowledgements}
The work of Hanna Wojew\'odka-\'Sci\k{a}\.zko has been partly supported by the National Science Centre of Poland, grant number 2018/02/X/ST1/01518.

\bibliography{references}
\bibliographystyle{plain}
\end{document}